\documentclass[12pt,reqno]{amsart}
\usepackage{fullpage}
\usepackage{amsfonts}
\usepackage{amssymb}
\usepackage{times}
\usepackage{graphicx}
\newtheorem{thm}{Theorem}[section]
\newtheorem{cor}[thm]{Corollary}

\begin{document}

\title[SUBCLASS OF HARMONIC UNIVALENT
FUNCTIONS ASSOCIATED WITH THE...]{SUBCLASS OF HARMONIC UNIVALENT
FUNCTIONS ASSOCIATED WITH THE GENERALIZED MITTAG-LEFFLER TYPE FUNCTIONS }
% or if you want, simply \title{Title of the article}
\author{Adnan Ghazy Alamoush \\
}
\date{\sl \small Faculty of Science, Taibah University, Saudi Aarabia.\\
}
\maketitle

\maketitle

\begin{abstract}
In the present paper, we introduce a new subclass of harmonic functions in the unit disc $U$ defined by using the generalized Mittag-Leffler type functions. Coefficient conditions, extreme points, distortion bounds, convex combination are studied.
\end{abstract}

\bigskip
{\it Keywords}: Harmonic functions, hypergeometric functions, Mittag-Leffler type functions.

\bigskip

AMS Mathematics Subject Classification: 30C45

\section{Introduction}
A continuous complex-valued function $f = u + iv $ defined in a simply complex domain $D$ is said to be harmonic in $D$. In any simply connected domain, we can write $f=h+\bar{g}$, where $h$ and $g$ are analytic in $D$. A necessary and sufficient condition for $f$ to be locally univalent and sense preserving in  $D$ is that$|h^{'}(z)|>|g^{'}(z)|,z\in D$.

\medskip\noindent
Clunie and Sheil-Small \cite{Clunie} introduced a class $S{H}$  of complex valued harmonic maps f which are univalent and sense-preserving in the open unit disk $U=\{{z: |z|<1}\}$ and assume a normalized representation  $f=h+\bar{g}$ where $f(0)=f_{z}(0)-1=0$.  Then for $f = h +\bar{g}\in S{H}$  we may express the analytic functions $ h$ and $g$ as

\begin{equation}\label{eg:}
        h(z)=z+\sum^{\infty}_{k=2} a_{k}z^{k},\    g(z)=\sum^{\infty}_{k=1} b_{k}z^{k}, |b_{1}|<1.
\end{equation}

\medskip\noindent
In 1984, Clunie\cite{Clunie} and Sheil-Small \cite{Sheil} investigated the class $S{H}$ as well as its geometric subclasses and obtained some coefficient bounds. Since then, there have been several related papers on $S{H}$ and its subclasses.
\medskip\noindent
Connectivity of geometric functions and hypergeometric functions with harmonic functions is seen through some of the these  papers [\cite{Avci}\cite{Adnan}, \cite{Ahuja}, \cite{Ahuja1},\cite{Ahuja2}]. The Mittag-Leffler and generalized Mittag-Leffler type functions was first introduced by the
swedish mathematician G. M. Mittag-Leffler \cite{Mittag} and also studied by Wiman \cite{Wiman}. It is a special function of $z\in C$ which depends on the complex parameter $\alpha$ and is defined by the power series
\begin{equation}\label{eg:}
      E_{\alpha}(z)=\sum^{\infty}_{k=0}\frac{z^{k}}{\Gamma(\alpha k+1)}, \alpha \in C, R(\alpha)>0, z\in C.
\end{equation}
A first generalization of $E_{\alpha}(z)$ introduced by Wiman \cite{Wiman}, is the two-parametric M-L function of $z \in C$, defined by the series
\begin{equation}\label{eg:}
      E_{\alpha,\beta}(z)=\sum^{\infty}_{k=0}\frac{z^{k}}{\Gamma(\beta+\alpha k)}, \alpha,\beta\in C, R(\alpha)>0, R(\beta)>0, z\in C.
\end{equation}
1971 Prabhakar \cite{Prabhakar} introduced a three-parametric generalization of $ \psi^{\gamma}_{\alpha,\beta}(z)$ defined in (3) as a kernel of certain fractional differential equations in terms of the series
\begin{equation}\label{eg:}
     E^{\gamma}_{\alpha,\beta}(z)=\sum^{\infty}_{k=0}\frac{(\gamma)_{k}z^{k}}{k!\Gamma (\beta+\alpha k)}, \alpha,\beta\in C, R(\alpha)>0, R(\beta)>0, z\in C.
\end{equation}
Due to its integral representation  $ E^{\gamma}_{\alpha,\beta}(z)$ is considered as a special case of Fox’s H-function as well as of Wright's generalized hypergeometric $_{p}\Psi_{q}$, so called
Fox-Wright psi function of $z\in C$.\
Later, Salim \cite{Tariq} introduced  the function in the form $ \psi^{\gamma}_{\alpha,\beta}(z)$ in the following form
\begin{equation}\label{eg:}
     E^{\gamma, \delta}_{\alpha,\beta}(z)=\sum^{\infty}_{k=0}\frac{(\gamma)_{k}z^{k}}{\Gamma (\beta+\alpha k)(\delta)_{k}},
\end{equation}
where $\alpha,\beta, \gamma, \delta\in C, \min(R(\alpha), R(\beta)>0, R(\gamma), R(\delta)>0), z\in C.$
Recently, Salim and Faraj \cite{Tariq1}  introduced a new generalization of Mittag-Leffler function associated with Weyl fractional integral and differential operators as follow
\begin{equation}\label{eg:}
     E^{\gamma, \delta,q}_{\alpha,\beta, p}(z)=\sum^{\infty}_{k=0}\frac{(\gamma)_{qk}z^{k}}{\Gamma (\beta+\alpha k)(\delta)_{pk}},
\end{equation}
where $ \alpha,\beta, \gamma, \delta\in C, \min(R(\alpha), R(\beta)>0, R(\gamma), R(\delta)>0), z\in C$, with $q,p\in \mathbb{R}_{+}$, $q\leq\Re(\alpha)+p$, and $(\gamma)_{pn}$ denotes an extended variant of the
Pochhammer symbol, defined by $(\gamma)_{qn}=\Gamma(\gamma+qn)/\Gamma(\gamma)$.

\medskip\noindent
Corresponding to $ E^{\gamma, \delta, q}_{\alpha,\beta, p}(z)$, we define the function $ \Theta^{\gamma, \delta, q}_{\alpha,\beta, p}(z)$ by
$$ \Theta^{\gamma, \delta, q}_{\alpha,\beta, p}(z)=z* E^{\gamma, \delta, q}_{\alpha,\beta, p}(z)$$
$$=z+\sum^{\infty}_{k=2}\frac{(\gamma)_{q(k-1)}}{\Gamma [\beta+\alpha (k-1)](\delta)_{p(k-1)}}z^{k}.$$
Now, for $f \in A, m\in \mathbb{N}$, we define the following differential operator: $\Phi^{m}_{\gamma, \delta, q,\alpha,\beta, p})f :A �\rightarrow A$ by
\[\Phi^{0}_{\gamma, \delta, q,\alpha,\beta, p}f(z)=f(z)*\Theta^{\gamma, \delta, q}_{\alpha,\beta, p}(z)\]
\[\Phi^{1}_{\gamma, \delta, q,\alpha,\beta, p}f(z)=z[f(z)*\Theta^{\gamma, \delta, q}_{\alpha,\beta, p}(z)]^{'}\]
\[.\]
\[.\]
\[.\]
\[(m+1)\Phi^{m+1}_{\gamma, \delta, q,\alpha,\beta, p}f(z)=z[\Phi^{m}_{\gamma, \delta, q,\alpha,\beta, p}]+m\Phi^{m}_{\gamma, \delta, q,\alpha,\beta, p},\ z\in U.\]
Thus it is obvious to see from above that
\begin{equation}\label{eg:}
\Phi^{m}_{\gamma, \delta, q,\alpha,\beta, p}=z+ \sum^{\infty}_{k=2}\frac{(m+1)_{k-1}}{(k-1)!}\frac{(\gamma)_{q(k-1)}}{\Gamma [\beta+\alpha (k-1)](\delta)_{p(k-1)}}a_{k}z^{k}.
\end{equation}

Note that, when $\alpha=0, \beta=\gamma=\delta=1$, we get Ruscheweyh Operator \cite{Ruscheweyh}.\\

Involving the generalized Mittag-Leffler function as defined in (6), for $0\leq\eta<1, m\in \mathbb{N}, n\in \mathbb{N}_{0}, m>n$ and $z\in U$,  we let $SH(m,n,\eta)$ denote the family of harmonic
functions f of the form (1) such that
\begin{equation}\label{eg:}
\Re\left(\frac{\Phi^{m}_{\gamma, \delta, q,\alpha,\beta, p}}{\Phi^{n}_{\gamma, \delta, q,\alpha,\beta, p}}\right)>\eta,
\end{equation}
where $\Phi^{m}_{\gamma, \delta, q,\alpha,\beta, p}=\Phi^{m}_{\gamma, \delta, q,\alpha,\beta, p}h(z)+(-1)^{m}\overline{\Phi^{m}_{\gamma, \delta, q,\alpha,\beta, p}}g(z).$\

We let the subclass $\overline{S}H(m,n,\eta)$ consist of harmonic functions $f_{m}=h+ \overline{g}_{m}$ in $\overline{S}H(m,n,\eta)$ so that $h$
and $g_{m}$ are of the form
\medskip\noindent
\begin{equation}\label{eg:}
  h(z)=z-\sum^{\infty}_{k=2}a_{k}z^{k},\ g_{m}(z)=(-1)^{m-1}\sum^{\infty}_{k=1}b_{k}z^{k}.
\end{equation}
The class $\overline{S}H(m,n,\eta)$ includes a variety of well-known subclasses of $SH$.
\medskip\noindent
The object of this paper is to examine some generalized Mittag-Leffler function inequalities
as a necessary and sufficient conditions for univalent harmonic analytic functions associated with
certain generalized Mittag-Leffler function to be in the function class $SH(m,n,\eta)$.
The coefficient condition for the function class $SH(m,n,\eta)$ is given. Furthermore, we determine extreme points, a distortion theorem, convolution conditions and convex combinations for the functions $f$ in $\overline{S}H(m,n,\eta)$.

\section{Coefficient bound}
We begin with a sufficient coefficient condition for functions $f$ in $SH(m,n,\eta)$.

\begin{thm}
Let $f = h + \overline{g}$ be given by (1). If
\begin{equation}\label{eg:}
\sum^{\infty}_{k=1}\left(\frac{\Phi^{m}_{\gamma, \delta, q,\alpha,\beta, p}-\eta \Phi^{n}_{\gamma, \delta, q,\alpha,\beta, p}}{1-\eta}|a_{k}|+\frac{\Phi^{m}_{\gamma, \delta, q,\alpha,\beta, p}-(-1)^{m-n}\eta \Phi^{n}_{\gamma, \delta, q,\alpha,\beta, p}}{1-\eta}|b_{k}|\right)\leq2.
\end{equation}
\end{thm}

\medskip\noindent
\textit{Proof}. If $z_{1}\neq z_{2}$ , then

\begin{center}
$\left|\frac{f(z_{1})-f(z_{2})}{h(z_{1})-h(z_{2})}\right|\geq1-\left|\frac{g(z_{1})-g(z_{2})}{h(z_{1})-h(z_{2})}\right|$
\end{center}
\begin{center}
$=1-\left|\frac{\sum^{\infty}_{k=1}b_{k}(z^{k}_{1}-z^{k}_{2})}{(z_{1}-z_{2})+
\sum^{\infty}_{k=2}a_{k}(z^{k}_{1}-z^{k}_{2})}\right|>1-
\frac{\sum^{\infty}_{k=1}k|b_{k}|}{1-\sum^{\infty}_{k=2}k|a_{k}|}$
\end{center}
\begin{center}
$\geq 1-\frac{\frac{\sum^{\infty}_{k=1}\Phi^{m}_{\gamma, \delta, q,\alpha,\beta, p}-(-1)^{m-n} \eta\Phi^{n}_{\gamma, \delta, q,\alpha,\beta, p}}{1-\eta}|b_{k}|}
{1-\sum^{\infty}_{k=2}\frac{\Phi^{m}_{\gamma, \delta, q,\alpha,\beta, p}-\eta\Phi^{n}_{\gamma, \delta, q,\alpha,\beta, p}}{1-\eta}|a_{k}|}\geq0,$
\end{center}
which proves univalence. Note that $f$ is sense-preserving in $U$. This is because

\begin{center}
$|h^{'}(z)|=1-\sum^{\infty}_{k=2}k|a_{k}| |z|^{k-1}$
\end{center}

\begin{center}
$>1-\sum^{\infty}_{k=2}\frac{\Phi^{m}_{\gamma, \delta, q,\alpha,\beta, p}-\eta \Phi^{n}_{\gamma, \delta, q,\alpha,\beta, p}}{1-\eta}|a_{k}|$
\end{center}

\begin{center}
$\geq\sum^{\infty}_{k=1}\frac{\Phi^{m}_{\gamma, \delta, q,\alpha,\beta, p}-(-1)^{m-n} \eta \Phi^{n}_{\gamma, \delta, q,\alpha,\beta, p}}{1-\eta}|b_{k}|$\end{center}

\begin{center}
$>\sum^{\infty}_{k=1}\frac{\Phi^{m}_{\gamma, \delta, q,\alpha,\beta, p}-(-1)^{m-n} \eta \Phi^{n}_{\gamma, \delta, q,\alpha,\beta, p}}{1-\eta}|b_{k}||z^{k}|$
\end{center}
\begin{center}
\begin{equation}\label{eg:}
\geq\sum^{\infty}_{k=1}k|b_{k}| |z|^{k-1}\geq|g'(z)|.
\end{equation}
\end{center}
Using the fact that $\Re(w) >\eta$ if and only if $ |1 -\eta + w| \geq |1 +\eta - w|$, it suffices to show that
\\

$\left|(1-\eta)\Phi^{n}_{\gamma, \delta, q,\alpha,\beta, p}(z)+\Phi^{m}_{\gamma, \delta, q,\alpha,\beta, p}(z)\right|$
\begin{equation}\label{eg:}
    -\left|(1+\eta)\Phi^{n}_{\gamma, \delta, q,\alpha,\beta, p}(z)-\Phi^{m}_{\gamma, \delta, q,\alpha,\beta, p}(z)
\right|>0.
\end{equation}
Substituting the value of $\Phi^{m}_{\gamma, \delta, q,\alpha,\beta, p}(z)$ and $\Phi^{n}_{\gamma, \delta, q,\alpha,\beta, p}(z)$  in (11) yields, by (9), that
\\
$$\left|(1-\eta)\Phi^{n}_{\gamma, \delta, q,\alpha,\beta, p}+\Phi^{m}_{\gamma, \delta, q,\alpha,\beta, p}\right|
    -\left|(1+\eta)\Phi^{n}_{\gamma, \delta, q,\alpha,\beta, p}-\Phi^{m}_{\gamma, \delta, q,\alpha,\beta, p}
\right|>0.$$

$=|(2-\eta)z+\sum^{\infty}_{k=2}[\Phi^{n}_{\gamma, \delta, q,\alpha,\beta, p}(1-\eta)+\Phi^{m}_{\gamma, \delta, q,\alpha,\beta, p}]a_{k}z^{k}$

\begin{center}
$+(-1)^{n}\sum^{\infty}_{k=1}[\Phi^{n}_{\gamma, \delta, q,\alpha,\beta, p}(1-\eta)+ (-1)^{m-n}\Phi^{m}_{\gamma, \delta, q,\alpha,\beta, p}]\overline{b_{k}z^{k}}$
\end{center}

$-|-\eta z+\sum^{\infty}_{k=2}[\Phi^{m}_{\gamma, \delta, q,\alpha,\beta, p}-(1+\eta)\Phi^{n}_{\gamma, \delta, q,\alpha,\beta, p}]a_{k}z^{k}|$

\begin{center}
$-(-1)^{n}\sum^{\infty}_{k=1}[(-1)^{m-n}\Phi^{m}_{\gamma, \delta, q,\alpha,\beta, p}-(1+\eta)\Phi^{n}_{\gamma, \delta, q,\alpha,\beta, p}]\overline{b_{k}z^{k}}|$
\end{center}
$\geq2(1-\eta)|z|-\sum^{\infty}_{k=2}2 [\Phi^{m}_{\gamma, \delta, q,\alpha,\beta, p}-\eta\Phi^{n}_{\gamma, \delta, q,\alpha,\beta, p}]|a_{k}|z|^{k}$
\begin{center}
$-\sum^{\infty}_{k=1}[(1+\eta)\Phi^{n}_{\gamma, \delta, q,\alpha,\beta, p}+(-1)^{m-n}\eta\Phi^{m}_{\gamma, \delta, q,\alpha,\beta, p}]|b_{k}||z|^{k}$
\end{center}
$-\sum^{\infty}_{k=1}|[(-1)^{m-n}\Phi^{m}_{\gamma, \delta, q,\alpha,\beta, p}-(1+\eta)\Phi^{n}_{\gamma, \delta, q,\alpha,\beta, p}]|b_{k}||z|^{k}$
$$=2(1-\eta)|z|$$
$$\left\{1-\sum^{\infty}_{k=2}\frac{\Phi^{m}_{\gamma, \delta, q,\alpha,\beta, p}-\eta \Phi^{n}_{\gamma, \delta, q,\alpha,\beta, p}}{1-\eta}|a_{k}||z|^{k-1}-\sum^{\infty}_{k=1}\frac{\Phi^{m}_{\gamma, \delta, q,\alpha,\beta, p}-(-1)^{m-n}\eta \Phi^{n}_{\gamma, \delta, q,\alpha,\beta, p}}{1-\eta}|b_{k}||z|^{k-1}\right\}$$

$$>2(1-\eta)|z|\left\{1-\left(\sum^{\infty}_{k=2}\frac{\Phi^{m}_{\gamma, \delta, q,\alpha,\beta, p}-\eta \Phi^{n}_{\gamma, \delta, q,\alpha,\beta, p}}{1-\eta}|a_{k}|+\sum^{\infty}_{k=1}\frac{\Phi^{m}_{\gamma, \delta, q,\alpha,\beta, p}-(-1)^{m-n}\eta \Phi^{n}_{\gamma, \delta, q,\alpha,\beta, p}}{1-\eta}|b_{k}|\right)\right\}.$$
This last expressions is non-negative by (10), and so the proof is complete.

\medskip\noindent
The harmonic function

\begin{equation}\label{eg:}f(z)=z+\sum^{\infty}_{k=2}\frac{1-\eta}{\Phi^{m}_{\gamma, \delta, q,\alpha,\beta, p}-\eta \Phi^{n}_{\gamma, \delta, q,\alpha,\beta, p}}x_{k}z^{k}+\sum^{\infty}_{k=1}\frac{1-\eta}{\Phi^{m}_{\gamma, \delta, q,\alpha,\beta, p}-(-1)^{m-n}\eta \Phi^{n}_{\gamma, \delta, q,\alpha,\beta, p}}\overline{y_{k}z^{k}},
\end{equation}
 where $m\in \mathbb{N}, n\in \mathbb{N}_{0}$ and $\sum^{\infty}_{k=2}|x_{k}|+\sum^{\infty}_{k=1}|y_{k}|=1$ shows that the coefficient bound given by (10) is sharp. The functions $f$ of the form (13) are $f\in SH(m,n,\eta)$, because
 $$\sum^{\infty}_{k=1}\left(\frac{\Phi^{m}_{\gamma, \delta, q,\alpha,\beta, p}-\eta \Phi^{n}_{\gamma, \delta, q,\alpha,\beta, p}}{1-\eta}|a_{k}|+\frac{\Phi^{m}_{\gamma, \delta, q,\alpha,\beta, p}-(-1)^{m-n}\eta \Phi^{n}_{\gamma, \delta, q,\alpha,\beta,p}}{1-\eta}|b_{k}|\right)=1+\sum^{\infty}_{k=2}|x_{k}|+\sum^{\infty}_{k=1}|y_{k}|=2.
$$

\medskip\noindent
In the following theorem, it is shown that the condition (10) is also necessary for functions $f_{m} = h + \overline{g_{m}}$,  where $h$ and $g_{m}$  are of the form (9).

\begin{thm}
Let $f_{m} = h + \overline{g_{m}}$ be given by (9). Then $f_{m}\in \overline{S}H(m,n, \eta)$ If and only if

\begin{equation}\label{eg}
\sum^{\infty}_{k=1}\left[\left(\Phi^{m}_{\gamma, \delta, q,\alpha,\beta, p}-\eta \Phi^{n}_{\gamma, \delta, q,\alpha,\beta, p}\right)|a_{k}|+\left(\Phi^{m}_{\gamma, \delta, q,\alpha,\beta, p}-(-1)^{m-n}\eta \Phi^{n}_{\gamma, \delta, q,\alpha,\beta, p}\right)|b_{k}|\right]\leq2(1-\eta).
\end{equation}
\end{thm}

\medskip\noindent
\textit{Proof}. Since $ \overline{S}H(m,n,\eta)\subset \overline{S}H(m,n,\eta)$, we only need to prove the "only if " part of Theorem 2.2. To this end, for functions  $f$ of the form (9), we notice that the condition (8) is equivalent to
\begin{center}
$  \Re\left\{\frac{(1-\eta)z-\sum_{k=2}^{\infty}\left(\Phi^{m}_{\gamma, \delta, q,\alpha,\beta, p}-\eta \Phi^{n}_{\gamma, \delta, q,\alpha,\beta, p}\right)a_{k}z^{k}+(-1)^{2m-1}\sum_{k=1}^{\infty}\left(\Phi^{m}_{\gamma, \delta, q,\alpha,\beta, p}-(-1)^{m-n}\eta \Phi^{n}_{\gamma, \delta, q,\alpha,\beta, p}\right)b_{k}\overline{z}^{k}}{z-\sum_{k=2}^{\infty}\Phi^{n}_{\gamma, \delta, q,\alpha,\beta, p}a_{k}z^{k}+(-1)^{m+n-1}\sum_{k=1}^{\infty}\Phi^{n}_{\gamma, \delta, q,\alpha,\beta, p}b_{k}\overline{z}^{k}}\right\}\geq0.$
\end{center}
\begin{equation}\label{eg:}
\end{equation}
The above condition (15) must hold for all values of $z$ on the positive real axes, where,
$0\leq|z|=\mu<1,$ we have

\begin{center}
$\frac{1-\eta-\sum_{k=2}^{\infty}\left(\Phi^{m}_{\gamma, \delta, q,\alpha,\beta, p}-\eta \Phi^{n}_{\gamma, \delta, q,\alpha,\beta, p}\right)a_{k}\mu^{k-1}-\sum_{k=1}^{\infty}\left(\Phi^{m}_{\gamma, \delta, q,\alpha,\beta, p}-(-1)^{m-n}\eta \Phi^{n}_{\gamma, \delta, q,\alpha,\beta, p}\right)b_{k}\mu^{k-1}}{1-\sum_{k=2}^{\infty}\Phi^{n}_{\gamma, \delta, q,\alpha,\beta, p}a_{k}\mu^{k-1}+(-1)^{m-n}\sum_{k=1}^{\infty}\Phi^{n}_{\gamma, \delta, q,\alpha,\beta, p}b_{k}\mu^{k-1}}\geq0.$
\end{center}
\begin{equation}\label{eg:}
\end{equation}

\medskip\noindent
If the condition (14) does not hold, then the numerator in (16) is negative for sufficiently close to 1. Hence there exists a $z_{0} =\mu_{0}  $ in $(0, 1)$ for which the quotient in (16) is negative. This contradicts the condition for $f_{j}\in \overline{S}H(m,n,\eta)$ and so the proof is complete.

\section{Distortion bounds}

In this section, we obtain distortion bounds for functions $f$ in $\overline{S}H(m,n,\eta).$

\begin{thm} Let $f_{j}\in \overline{S}H(m,n,\eta)$, then for $|z|<1$, we have

\begin{center}
$|f_{k}(z)|\geq(1-|b_{1}|)\mu-\frac{1}{\Upsilon_{2}}\left(\frac{(1-\eta)}{[\Upsilon_{2}]^{m-n}-\eta}
-\frac{1-(-1)^{m-n}\eta}{[\Upsilon_{2}]^{m-n}-\eta}|b_{1}|\right)\mu^{2}$
\end{center}
and

\begin{center}
$|f_{k}(z)|\geq(1-|b_{1}|)\mu-\frac{1}{\Upsilon_{2}}\left(\frac{(1-\eta)}{[\Upsilon_{2}]^{m-n}-\eta}
-\frac{1-(-1)^{m-n}\eta}{[\Upsilon_{2}]^{m-n}-\eta}|b_{1}|\right)\mu^{2},$
\end{center}
where $\Upsilon_{2}=(m+1)\left(\frac{(\gamma)_q}{\Gamma [\beta+\alpha ](\delta)_{p}}\right).$
\end{thm}

\medskip\noindent
\textit{Proof.} We only prove the left-hand inequality. The proof for the right-hand inequality
is similar and is thus omitted. Let  $f_{m}\in \overline{S}H_{\alpha, \beta, p,n}^{\gamma,\delta,q,m}(\eta)$.  Taking the absolute value of $f_{m}$, we have

\begin{center}
$|f_{m}(z)|=|z-\sum_{k=2}^{\infty}a_{k}z^{k}+(-1)^{k}\sum_{k=1}^{\infty}b_{k}\overline{z}^{k}|$
\end{center}
\begin{center}

$\leq(1+|b_{1}|)\mu+\sum_{k=2}^{\infty}(|a_{k}|+|b_{k}|)\mu^{n}$
\end{center}

\begin{center}
    $\leq(1+|b_{1}|)\mu+\sum_{k=2}^{\infty}(|a_{k}|+|b_{k}|)\mu^{2}$
\end{center}
\begin{center}
$=(1+|b_{1})\mu+\frac{1-\eta}{[\Upsilon_{2}]^{m}([\Upsilon_{2}]^{m-n}-\eta)}\sum^{\infty}_{k=2}(\frac{\Upsilon_{2}]^{m}
([\Upsilon_{2}]^{m-n}-\eta)}{1-\gamma}(|a_{k}|+|b_{k}|)\mu^{2}$
\end{center}

\begin{center}
$\leq(1+|b_{1})\mu+\frac{(1-\eta)\mu^{2}}{[\Upsilon_{2}]^{m}([\Upsilon_{2}]^{m-n}-\eta)}
\sum^{\infty}_{k=2}\left(\frac{\Phi^{m}_{\gamma, \delta, q,\alpha,\beta, p}-\eta \Phi^{n}_{\gamma, \delta, q,\alpha,\beta, p}}{1-\eta}a_{k}+\frac{\Phi^{m}_{\gamma, \delta, q,\alpha,\beta, p}-(-1)^{m-n}\eta \Phi^{n}_{\gamma, \delta, q,\alpha,\beta, p}}{1-\eta}b_{k}\right)$\end{center}

\begin{center}
$\leq(1+|b_{1})\mu+\frac{(1-\eta)\mu^{2}}{[\Upsilon_{2}]^{m}([\Upsilon_{2}]^{m-n}-\eta)}\left(\frac{\Phi^{m}_{\gamma, \delta, q,\alpha,\beta, p}-\eta \Phi^{n}_{\gamma, \delta, q,\alpha,\beta, p}}{1-\eta}a_{k}+\frac{\Phi^{m}_{\gamma, \delta, q,\alpha,\beta, p}-(-1)^{m-n}\eta \Phi^{n}_{\gamma, \delta, q,\alpha,\beta, p}}{1-\eta}b_{k}\right)\mu^{2}$.
\end{center}
The following covering result follows from the left hand inequality in Theorem 3.1.
\section{Convolution, convex combinations and extreme points}
In this section, we show the class $SH(m,n,\eta)$  is invariant under convolution and convex combination.

\medskip\noindent
For harmonic functions $f$ of the form

\begin{center}
    $f_{m}(z)=z-\sum_{k=2}^{\infty}|a_{k}|z^{k}+(-1)^{m-1}\sum_{k=1}^{\infty}|b_{k}|\bar{z}^{k}$
\end{center}
and
\begin{center}
    $F_{m}(z)=z-\sum_{k=2}^{\infty}|A_{k}|z^{k}+(-1)^{m-1}\sum_{k=1}^{\infty}|B_{k}|\bar{z}^{k},$
\end{center}
we define the convolution of $f_{m}(z)$ and $F_{m}(z)$ as

\begin{equation}\label{eg:}
    (f_{m}*F_{m})(z)=f_{m}(z)*F_{m}(z)
    =z-\sum_{k=2}^{\infty}|a_{k}||A_{k}|z^{n}+(-1)^{k}\sum_{k=1}^{\infty}|b_{k}||B_{k}|\bar{z}^{k}.
\end{equation}

\begin{thm}
For $0\leq\rho\leq\eta<1$, let $f_{m}\in \overline{S}H(m,n,\eta)$ and
$F_{m}\in \overline{S}(m,n,\rho)$. Then the convolution
\begin{center}
    $f_{m}*F_{m}\in \overline{S}H(m,n,\eta)\subset \overline{S}(m,n,\rho).$
\end{center}
\end{thm}

\medskip\noindent
\textit{Proof. }Then the convolution $f_{m}*F_{m}$ is given by (17). We wish to show that the
coefficients of $f_{m}*F_{m}$ satisfy the required condition given in Theorem 4.1.
For $F_{m}\in \overline{S}H(m,n,\rho)$, we note that $|A_{k}|\leq1$ and $|B_{k}|\leq1$. Now, for the convolution function $f_{m}*F_{m}$, we obtain
\begin{center}
$\sum^{\infty}_{k=2}\frac{\Phi^{m}_{\gamma, \delta, q,\alpha,\beta, p}-\rho \Phi^{n}_{\gamma, \delta, q,\alpha,\beta, p}}{1-\rho}a_{k}A_{k}+\sum^{\infty}_{k=1}\frac{\Phi^{m}_{\gamma, \delta, q,\alpha,\beta, p}-(-1)^{m-n}\rho \Phi^{n}_{\gamma, \delta, q,\alpha,\beta, p}}{1-\rho}b_{k}B_{k}$
\end{center}

\begin{center}
$\leq\sum^{\infty}_{k=2}\frac{\Phi^{m}_{\gamma, \delta, q,\alpha,\beta, p}-\rho \Phi^{n}_{\gamma, \delta, q,\alpha,\beta, p}}{1-\rho}a_{k}+\sum^{\infty}_{k=1}\frac{\Phi^{m}_{\gamma, \delta, q,\alpha,\beta, p}-(-1)^{m-n}\rho \Phi^{n}_{\gamma, \delta, q,\alpha,\beta, p}}{1-\rho}b_{k}$
\end{center}
\begin{center}
$\leq\sum^{\infty}_{k=2}\frac{\Phi^{m}_{\gamma, \delta, q,\alpha,\beta, p}-\mu \Phi^{n}_{\gamma, \delta, q,\alpha,\beta, p}}{1-\mu}a_{k}+\sum^{\infty}_{k=1}\frac{\Phi^{m}_{\gamma, \delta, q,\alpha,\beta, p}-(-1)^{m-n}\mu \Phi^{n}_{\gamma, \delta, q,\alpha,\beta, p}}{1-\mu}b_{k}\leq1.$
\end{center}

Since $0\leq\rho\leq\eta<1$, and $f_{m}\in \overline{S}H(m,n,\eta)$, then
$f_{m}*F_{m}\in \overline{S}H(m,n,\eta)\subset \overline{S}H(m,n,\rho)$.

\medskip\noindent
Next we discuss the convex combinations of the class $\overline{S}H(m,n,\eta)$.
\begin{thm}
The family $\overline{S}H(m,n,\eta)$ is closed under convex combination.
\end{thm}

\medskip\noindent
\textbf{Proof.} For $i = 1, 2, . . . $, suppose that $f_{m,j}\in \overline{S}H(m,n,\eta)$, where
\begin{equation}\label{eg;}
     f_{m,j}(z)=z-\sum_{k=2}^{\infty}|a_{k,j}|z^{k}+(-1)^{m-1}\sum_{k=1}^{\infty}|b_{k,j}|\bar{z}^{k}.
\end{equation}
Then by Theorem 2.1,
\begin{center}
$\sum^{\infty}_{k=1}\left(\frac{\Phi^{m}_{\gamma, \delta, q,\alpha,\beta, p}-\eta \Phi^{n}_{\gamma, \delta, q,\alpha,\beta, p}}{1-\eta}a_{k,j}+\frac{\Phi^{m}_{\gamma, \delta, q,\alpha,\beta, p}-(-1)^{m-n}\eta \Phi^{n}_{\gamma, \delta, q,\alpha,\beta, p}}{1-\eta}b_{k,j}\right) \leq2.$
\end{center}

\begin{equation}\label{eg:}
\end{equation}
For $\sum^{\infty}_{j=1}t_{j}=1$, $0\leq t_{j}<1 $, the convex combination of $f_{m,j}$ may be written as

\begin{center}
    $\sum^{\infty}_{j=1}t_{j}f_{m,j}(z)=z-\sum^{\infty}_{k=2}\left(\sum^{\infty}_{j=1}t_{j}|a_{k,j}|\right)z^{k}+
    (-1)^{m-1}\sum^{\infty}_{k=1}\left(\sum^{\infty}_{j=1}t_{j}|b_{k,j}|\right)\bar{z}^{k}.$
\end{center}
Then by (19),
\begin{center}
$\sum^{\infty}_{k=1}\left(\frac{\Phi^{m}_{\gamma, \delta, q,\alpha,\beta, p}-\eta \Phi^{n}_{\gamma, \delta, q,\alpha,\beta, p}}{1-\eta}\sum^{\infty}_{j=1}t_{j}a_{k,j}+
\frac{\Phi^{m}_{\gamma, \delta, q,\alpha,\beta, p}-(-1)^{m-n}\eta \Phi^{n}_{\gamma, \delta, q,\alpha,\beta, p}}{1-\eta}\sum^{\infty}_{j=1}t_{j}b_{k,j}\right)$
\end{center}
\begin{center}
    $=\sum^{\infty}_{i=1}t_{i}\left\{\sum^{\infty}_{k=1}\left(\frac{\Phi^{m}_{\gamma, \delta, q,\alpha,\beta, p}-\eta \Phi^{n}_{\gamma, \delta, q,\alpha,\beta, p}}{1-\eta}a_{k,j}+\frac{\Phi^{m}_{\gamma, \delta, q,\alpha,\beta, p}-(-1)^{m-n}\eta \Phi^{n}_{\gamma, \delta, q,\alpha,\beta, p}}{1-\eta}b_{k,j}\right) \right\}$

\end{center}

\begin{center}
    $\leq2\sum^{\infty}_{j=1}t_{j}=2$,
\end{center}
 and therefore
\begin{center}
    $\sum^{\infty}_{j=1}t_{j}f_{m,j}(z)\in \overline{S}H(m,n,\eta).$
\end{center}
\begin{cor}
The class $\overline{S}H(m,n,\eta)$ is closed under convex linear combinations.
\end{cor}

\medskip\noindent
\textit{Proof.} Let the functions $f_{m,j} (z)\ (j = 1, 2 . . . ,m)$ defined by (18) be in the class $\overline{S}H(m,n,\eta).$ Then the function $\varpi(z)$ defined by
\begin{equation}\label{eg:}
    \varpi(z)=\mu f_{m,j}(z)+(1-\mu)f_{m,j}(z), 0\leq\mu\leq1,
\end{equation}
is in the class $\overline{S}H(m,n,\eta)$.

\medskip\noindent
Next we determine the extreme points of closed convex hulls of $\overline{S}H(m,n,\eta)$, denoted by clco\ $\overline{S}H(m,n,\eta).$

\begin{thm}
Let $f_{m}$ be given by (7). Then $f_{m}\in \overline{S}H(m,n,\eta)$ if and only if
\begin{equation}\label{eg:}
    f_{m}(z)=\sum^{\infty}_{k=1}\left[X_{k}h_{k}(z)+Y_{k}g_{m_{k}}(z)\right],
\end{equation}
where

\end{thm}
\begin{center}
    $h_{1}(z)=1,\ h_{k}(z)=z-\frac{1-\eta}{\Phi^{m}_{\gamma, \delta, q,\alpha,\beta, p}-\eta \Phi^{n}_{\gamma, \delta, q,\alpha,\beta, p}}z^{k},  (k=2,...),$
\end{center}
\begin{center}
   $ g_{m_{k}}(z)=z+(-1)^{m-1}\frac{1-\eta}{\Phi^{m}_{\gamma, \delta, q,\alpha,\beta, p}-(-1)^{m-n}\eta \Phi^{n}_{\gamma, \delta, q,\alpha,\beta, p}}\bar{z}^{n},  (k=1,2,...),$
\end{center}

$X_{1}=1-\sum^{\infty}_{k=2}(X_{k}+Y_{k})\geq0,\ X_{k}\geq0,\ Y_{k}\geq0$. In particular, the extreme points of $\overline{S}H(m,n,\eta)$
are $\left\{h_{k}\right\} \ and\  \left\{g_{m_{k}}\right\}$.

\medskip\noindent
\textit{Proof.} For functions $f_{m}$ of the form (21) we have

\begin{center}
    $    f_{m}(z)=\sum^{\infty}_{k=1}\left[X_{k}h_{k}(z)+Y_{k}g_{m_{k}}(z)\right]$
\end{center}
\begin{center}
    $=\sum^{\infty}_{k=1}(X_{k}+Y_{k})z-\sum^{\infty}_{k=2}
        \frac{1-\eta}{\Phi^{m}_{\gamma, \delta, q,\alpha,\beta, p}-\eta \Phi^{n}_{\gamma, \delta, q,\alpha,\beta, p}}X_{k}z^{k}$
\end{center}
\begin{center}
    $+(-1)^{m-1}\sum^{\infty}_{n=1}\frac{1-\eta}{\Phi^{m}_{\gamma, \delta, q,\alpha,\beta, p}-(-1)^{m-n}\eta \Phi^{n}_{\gamma, \delta, q,\alpha,\beta, p}}Y_{k}\bar{z}^{k}.$
\end{center}
Then
\begin{center}
    $\sum^{\infty}_{k=2}\left(
        \frac{\Phi^{m}_{\gamma, \delta, q,\alpha,\beta, p}-\eta \Phi^{n}_{\gamma, \delta, q,\alpha,\beta, p}}{1-\eta}\right)
        \frac{1-\eta}{\Phi^{m}_{\gamma, \delta, q,\alpha,\beta, p}-\eta \Phi^{n}_{\gamma, \delta, q,\alpha,\beta, p}}X_{k}$
\end{center}
\begin{center}
    $\sum^{\infty}_{k=1}\left(\frac{\Phi^{m}_{\gamma, \delta, q,\alpha,\beta, p}-(-1)^{m-n}\eta \Phi^{n}_{\gamma, \delta, q,\alpha,\beta, p}}{1-\eta}\right)\frac{1-\eta}{\Phi^{m}_{\gamma, \delta, q,\alpha,\beta, p}-(-1)^{m-n}\eta \Phi^{n}_{\gamma, \delta, q,\alpha,\beta, p}}Y_{k}.$
\end{center}\begin{center}
    $=\sum^{\infty}_{n=2}X_{n}+\sum^{\infty}_{n=1}Y_{n}=1-X_{1}\leq1,$
\end{center}
\begin{equation}\label{eg:}
\end{equation}
and so $f_{m}\in clco\ \overline{S}H(m,n,\eta).$

\medskip\noindent
Conversely, suppose that $f_{m}\in clco\ \overline{S}H(m,n,\eta).$ Setting

\begin{center}
    $X_{k}=\frac{\Phi^{m}_{\gamma, \delta, q,\alpha,\beta, p}-\eta \Phi^{n}_{\gamma, \delta, q,\alpha,\beta, p}}{1-\eta}a_{k},\ 0\leq X_{k}\leq 1, k=2,...\ ,$
\end{center}

\begin{center}
    $Y_{k}=\frac{\Phi^{m}_{\gamma, \delta, q,\alpha,\beta, p}-(-1)^{m-n}\eta \Phi^{n}_{\gamma, \delta, q,\alpha,\beta, p}}{1-\eta}b_{k},\ 0\leq Y_{k}\leq 1, k=1,2,... \ ,$
\end{center}
and $X_{1}=1-\sum^{\infty}_{k=2}X_{k}+\sum^{\infty}_{k=1}Y_{k}$ , and note that, by Theorem 2.1, $x_{1} \geq 0$. Consequently, we obtain
\

\begin{center}
    $f_{m}(z)=\sum^{\infty}_{k=1}\left(h_{k}(z)X_{k}+g_{m_{k}}(z)Y_{k}\right),$ as required.
\end{center}
\

Using Corollary 4.3 we have $ clco\ \overline{S}H(m,n,\eta)=\overline{G}_{S}H(m,n,\eta)$. Then the statement of
Theorem 4.4 is true for $f\in \overline{G}_{S}H(m,n,\rho).$

\ \\

    $Adnan\ Ghazy\ Alamoush $

 \small Faculty of Science, Taibah University, Saudi Aarabia.\\

     email: {\em adnan--omoush@yahoo.com}
%
%-------------------------------------------------------------------
%
\end{document}